\documentclass[12pt,a4paper]{article}

\usepackage{amsfonts}

\renewcommand{\theequation}{\arabic{section}.\arabic{equation}}

\def\G{{\cal G}}                        %
\def\A{{\cal A}}                        %
\def\K{{\cal K}}                        %
\def\E{{\cal E}}                        %
\def\R{{\cal R}}                        %
\def\B{{\cal B}}                        %
\def\cK{{\check{\cal K}}}               %
\def\ad{{\mathrm{ad}\,}}                %
\def\bZ{{\mathbb Z}}                    %
\def\cZ{{\cal Z}}                       %
\def\bC{{\mathbb C}}                    %
\def\bR{{\mathbb R}}                    %
\def\bN{{\mathbb N}}                    %
\def\T{{\cal T}}                        %

\begin{document}

\begin{titlepage} 
\vspace*{0.5cm}
\begin{center}

{\large\bf  Generalizations of Felder's elliptic dynamical $r$-matrices  \\

associated with twisted  loop algebras of self-dual Lie algebras}

\end{center}

\vspace{1.5cm}

\begin{center}
{\large L. Feh\'er\footnote{Corresponding author, E-mail: lfeher@sol.cc.u-szeged.hu }
and B.G. Pusztai}
\end{center}
\bigskip
\begin{center}
Department of Theoretical Physics,  University of Szeged \\
Tisza Lajos krt 84-86, H-6720 Szeged, Hungary \\
\end{center}

\vspace{2.2cm}

\begin{abstract}  

A dynamical $r$-matrix is associated with every self-dual Lie algebra $\A$ which is
graded by finite-dimensional subspaces as $\A=\oplus_{n \in \cZ} \A_n$, where 
$\A_n$ is dual to $\A_{-n}$ with respect to the
invariant scalar product on $\A$, and $\A_0$ admits a nonempty open subset 
$\check \A_0$ for which $\ad \kappa$ is 
invertible on $\A_n$  if $n\neq 0$ and $\kappa \in \check \A_0$.
Examples are furnished by taking $\A$ to be an affine Lie algebra 
obtained from the central extension of a twisted loop algebra
$\ell(\G,\mu)$ of a finite-dimensional self-dual Lie algebra $\G$.
These $r$-matrices, $R: \check \A_0 \rightarrow \mathrm{End}(\A)$,
yield generalizations of the basic trigonometric dynamical $r$-matrices that,
according to Etingof and Varchenko, are associated with the Coxeter automorphisms 
of the simple Lie algebras, and are related to Felder's elliptic
$r$-matrices by evaluation homomorphisms of $\ell(\G,\mu)$ into $\G$.
The spectral-parameter-dependent  
dynamical $r$-matrix that corresponds analogously 
to an arbitrary scalar-product-preserving finite order automorphism of a
self-dual Lie algebra is here calculated explicitly.

\end{abstract}

\end{titlepage}

\section{Introduction}
\setcounter{equation}{0}

The classical dynamical Yang-Baxter equation (CDYBE) introduced 
in its general form by Etingof and Varchenko \cite{EV} is a
remarkable generalization of the CYBE.
Currently we are witnessing intense research on  
the theory and the applications of the CDYBE to integrable systems 
\cite{Feld,ABB,Xu}.  For a review, see \cite{ES1}.

The aim of this paper is to study infinite-dimensional
generalizations of a certain class of finite-dimensional classical dynamical 
$r$-matrices.
Next we briefly recall these finite-dimensional $r$-matrices, which  
appear naturally in the chiral WZNW model (see e.g. \cite{Dub}
and references therein).
  
Let $\A$ be a finite-dimensional complex Lie algebra equipped with a 
nondegenerate, symmetric, 
invariant bilinear form $\langle\ ,\ \rangle$.  
Such a Lie algebra is called self-dual \cite{Figu}.
Consider a self-dual subalgebra $\K\subset \A$, on which 
$\langle\ ,\ \rangle$ remains
nondegenerate.  Introduce the complex analytic functions $f$ and $F$ by 
\begin{equation}
f: z \mapsto \frac{1}{2}\coth \frac{z}{2}- \frac{1}{z},
\qquad
F: z \mapsto \frac{1}{2}\coth \frac{z}{2}.
\label{1.1}\end{equation}
Suppose that $\cK$ is a nonempty open subset of $\K$ on 
which the operator valued function $R: \cK \rightarrow \mathrm{End}(\A )$
is defined by 
\begin{equation}
R(\kappa):=\left\{
\begin{array}{cc} 
f(\ad \kappa) &\mbox{on $\K$}\\
F(\ad \kappa) 
 &\mbox{on $\K^\perp$} 
\end{array}\right. \qquad \forall \kappa\in \cK.
\label{1.2}\end{equation}
The decomposition $\A=\K + \K^\perp$ is induced by $\langle\ ,\ \rangle$.
$R(\kappa)$ is a well defined linear operator 
on $\A$ if and only if the spectrum of $\ad \kappa$, acting on $\A$,  
does not intersect $2\pi i \bZ^*$, and    
$\left(\ad \kappa\right)\vert_{\K^\perp}$ is 
invertible\footnote{The set of integers is denoted by $\bZ$, 
$\bZ^*=\bZ\setminus \{0\}$,  
and $\bN$ denotes the positive integers.}.
On $\cK\subset \K$ subject to these conditions, 
the following (modified) version of the CDYBE holds:
\begin{eqnarray}
&&[ R X, R Y] -R( [X, R Y]+ [RX,Y]) 
+\langle X, (\nabla R) Y\rangle + (\nabla_{Y_\K} R) X -
(\nabla_{X_\K} R) Y \nonumber\\
&&\qquad\qquad \qquad = -\frac{1}{4} [X,Y], \qquad \forall X, Y\in \A.
\label{1.3}\end{eqnarray}
Here the `dynamical variable' $\kappa$ is suppressed for brevity,
$\forall X\in \A$ is decomposed as $X=X_\K + X_{\K^\perp}$, 
and 
\begin{equation}
(\nabla_T R)(\kappa):= \frac{d}{dt} R(\kappa + t T) \vert_{t=0}
\qquad
\forall T\in \K,\quad \kappa\in \check \K,
\label{1.4}\end{equation} 
\begin{equation}
\langle X, (\nabla R)(\kappa) Y\rangle  := 
\sum_i K^i \langle X, (\nabla_{K_i} R)(\kappa) Y \rangle,
\qquad \forall X, Y\in \A,
\label{1.5}\end{equation}
where $ K_i$ and $K^i$  denote dual bases of $\K$,
$\langle K_i, K^j\rangle = \delta_i^j$.
$R(\kappa)$ is antisymmetric,
$\langle R(\kappa) X, Y\rangle = - \langle X, R(\kappa) Y\rangle$,
and is $\K$-equivariant in the sense that
\begin{equation}
(\nabla_{[T, \kappa]  } R)(\kappa) = [\ad T, R(\kappa)],
\qquad \forall T\in \K, \kappa\in \check \K.
\label{1.6}\end{equation}
These properties of $R$ have been established in this general 
setting in \cite{Dub,FGP1}.
In various special cases --- in particular the case $\K=\A$ --- they 
were proved earlier in \cite{EV,AM,BFP}.
If one introduces $r^\pm: \cK \rightarrow \A\otimes \A$ by
\begin{equation}
r^\pm(\kappa):= \left(R(\kappa) T_\alpha\right) \otimes T^\alpha 
\pm \frac{1}{2} T_\alpha \otimes T^\alpha,  
\label{1.7}\end{equation}
where $\{T_\alpha\}$ and $\{T^\alpha\}$ are dual bases of $\A$,
and uses the identification 
$\K\simeq \K^*$ induced by $\langle\ ,\ \rangle$,
then the above properties of $R$ become the CDYBE for $r^\pm$
with respect to the pair $\K\subset \A$ as defined in \cite{EV}.

It is natural to suspect that  whenever (\ref{1.2}) is a well defined
formula, the resulting $r$-matrix always satisfies (\ref{1.3}).
For this it is certainly not necessary to assume that $\A$ is finite dimensional.
For example, Etingof and Varchenko \cite{EV} verified the CDYBE in the situation
for which $\A$ is an affine Lie algebra based on a simple Lie algebra
and $\K\subset \A$ is a Cartan subalgebra. 
Moreover, by applying evaluation homomorphisms to these
$r$-matrices they recovered Felder's celebrated 
spectral-parameter-dependent elliptic dynamical $r$-matrices \cite{Feld}.
Without presenting proofs, this construction was generalized in \cite{FGP1}
to any affine Lie algebra, $\A(\G,\mu)$, defined by adding the
derivation to the central extension of a twisted loop algebra,
$\ell(\G,\mu)$, based on an appropriate automorphism, $\mu$, 
of a  self-dual Lie algebra, $\G$. 
Namely,  such an affine Lie algebra automatically comes equipped
with the integral gradation associated with the powers of the
loop parameter, and it can be shown that (\ref{1.2})  
provides a solution of (\ref{1.3}) if one takes $\K$ to be the
grade zero subalgebra in this gradation. 
In this paper,  this solution will arise as a special case 
of a general theorem, which ensures the validity of (\ref{1.3})
for (\ref{1.2}) under the assumption that $\K=\A_0$ where 
 $\A=\oplus_{n\in \cZ} \A_n$ is graded by
finite-dimensional subspaces and carries an invariant scalar product that
is compatible with the grading in the sense that $\A_n \perp \A_m$ 
unless $(n+m)=0$.
Here $\cZ$ is some abelian group, in our examples $\cZ=\bZ$.
The precise statement, 
which is our first main result, is given by theorem 1 in section 2.
We shall use this result to obtain dynamical $r$-matrices
on the twisted loop algebras $\ell(\G,\mu)$ with the dynamical
variable lying in the fixed point set $\G_0\subset \G$
of the automorphism $\mu$ of $\G$.
By means of evaluation homomorphisms, these 
$r$-matrices then yield spectral-parameter-dependent 
$\G\otimes \G$-valued dynamical $r$-matrices generalizing Felder's
elliptic $r$-matrices.
The latter are recovered if $\G$ is taken to be a simple Lie algebra
and $\mu$ a Coxeter automorphism, consistently with 
the derivation found in \cite{EV}.
The existence of the above-mentioned family of elliptic dynamical $r$-matrices 
was announced in \cite{FGP1}.
Our second main result is their derivation presented in section 3.
See in particular proposition 2 and proposition 3 in subsection 3.3.
We shall also find a relationship between  
the underlying $\ell(\G,\mu) \otimes \ell(\G,\mu)$-valued 
$r$-matrices with dynamical variables in $\G_0$, 
and certain $\G\otimes \G$-valued dynamical $r$-matrices on $\G_0$ 
introduced in \cite{ES2}.
This is contained in an appendix. 

Before turning to the main text, the reader may consult the concluding  
section, where the results are summarized once more and some comments
are offered on the possible applications of our dynamical $r$-matrices.

\section{$r$-matrices on graded, self-dual Lie algebras}
\setcounter{equation}{0}

In this section we apply formula (\ref{1.2}) to infinite-dimensional Lie algebras 
that are decomposed into finite-dimensional subspaces in such a way that 
the $r$-matrix
leaves these subspaces invariant. The definition of the $r$-matrix on 
these subspaces
will be given in terms of the well known 
holomorphic functional calculus of linear operators \cite{DS}.

Let $A\in \mathrm{End}(V)$ be a linear operator on a finite-dimensional 
complex vector 
space $V$. Denote by $\sigma_A$ the spectrum (set of eigenvalues) of $A$.
Consider a holomorphic complex function $H$ defined on some open domain 
containing $\sigma_A$,
and take  $\Gamma$ to be a contour that lies in this domain and encircles 
each eigenvalue 
of $A$ with orientation used in Cauchy's theorem. 
Then the operator $H(A)\in \mathrm{End}(V)$ may be defined by 
\begin{equation}
H(A) := \frac{1}{2\pi i} \oint_\Gamma dz H(z) (z I_V - A)^{-1},
\label{2.1}\end{equation}
where $I_V$ is the identity operator on $V$.  
This definition can be converted into an explicit formula by means of the 
Jordan decomposition of $A$, which shows that $H(A)$ only depends on the
derivatives $H^{(k)}(z_i)$ for $z_i\in \sigma_A$ up to a finite order.
For example, if $A v = z_i v$ then $H(A) v= H(z_i) v$ and 
for the constant function $H_c(z)\equiv c$ one obtains $H_c(A)= c I_V$. 
Furthermore, if the power series expansion 
$H(z)=\sum_{k=0}^\infty c_k z^k$ is valid 
in a neighbourhood of $\sigma_A$, then
$H(A)= \sum_{k=0}^\infty c_k A^k$.
An important rule of this functional calculus is that if $H_3=H_1 H_2$ 
on some admissible
domain, then $H_3(A)= H_1(A) H_2(A)$. 
See e.g.~chapter VII of \cite{DS}.

We now consider a complex Lie algebra  $\A$ equipped with a gradation 
based on some abelian group $\cZ$.
We use the additive notation to denote the group operation on $\cZ$. 
The zero as a number and the unit element of $\cZ$ are both 
denoted simply by $0$, but this should not lead to any confusion. 
We assume that as a linear space
\begin{equation}
\A= \oplus_{n\in \cZ} \A_n,
\quad
0\leq \mathrm{dim}(\A_n)<\infty,
\quad \mathrm{dim}(\A_0)\neq 0, 
\label{2.2}\end{equation}
and 
\begin{equation}
[\A_m, \A_n]\subset \A_{m+n}
\qquad
\forall m,n\in \cZ.
\label{2.3}\end{equation}
The elements of $\A$ are finite linear 
combinations of the elements of the homogeneous subspaces, and 
we permit the possibility that $\mathrm{dim}(\A_n)=0$ for some $n\in \cZ$.  
We further assume that $\A$ has a nondegenerate, symmetric,
invariant bilinear form
$\langle\ ,\ \rangle: \A \times \A \rightarrow \bC$, which is compatible 
with the gradation in the sense that
\begin{equation}
\A_m \perp \A_n 
\quad\hbox{unless}\quad (m+n)=0.
\label{2.4}\end{equation}
This means that if $(m+n)\neq 0$ then 
$\langle X, Y\rangle=0$ for any $X\in \A_m$, $Y\in \A_n$, and 
the dual space of 
$\A_n$ can be identified with $\A_{-n}$ by means of the pairing
given by $\langle\ ,\ \rangle$. 
In particular, $\A_0$ is a finite-dimensional self-dual subalgebra of $\A$.
Since $[\A_0, \A_n]\subset \A_n$ and $\A_n$ is finite dimensional,
$e^{\ad \kappa}$ is a well defined linear operator on $\A$ for any $\kappa\in \A_0$.
The invariance of the bilinear form,
$\langle [X,Y], Z \rangle + \langle Y, [X,Z] \rangle =0$,
$\forall X,Y,Z\in \A$,
implies that $\langle e^{\ad \kappa} Y, e^{\ad \kappa} Z\rangle = \langle Y, Z\rangle$ 
for any $Y,Z\in \A$ and $\kappa\in \A_0$.

Now we wish to apply formula (\ref{1.2}) to
\begin{equation}
\K:= \A_0,
\qquad
\K^\perp= \oplus_{n\in \cZ\setminus\{0\}} \A_n.
\label{2.5}\end{equation}
For any $\kappa\in \K$ and $n\in \cZ$, introduce 
$(\ad \kappa)_n:= \ad \kappa\vert_{\A_n}$ and let $\sigma_\kappa^n$ denote
the spectrum of this finite-dimensional linear operator
($\sigma_\kappa^n= \emptyset$ if $\mathrm{dim}(\A_n)=0$).
Our crucial assumption is that there exists a {\em nonempty, open}
subset $\check \K\subset \K$ for which 
\begin{equation}
\sigma_\kappa^n \cap 2\pi i \bZ =\emptyset 
\quad\forall n\neq 0\quad\hbox{and}\quad 
\sigma^0_\kappa \cap 2\pi i \bZ^* =\emptyset
\quad \forall \kappa \in \cK, 
\label{2.6}\end{equation}
where $\bZ$ and $\bZ^*$ are the set of all integers, 
and nonzero integers, respectively.
It is clear that if such a $\check \K$ exists, 
then there exists also a maximal one.  
If this assumption is satisfied, then we can define the map
$R: \check \K \rightarrow \mathrm{End}(\A)$ by requiring  that
the homogeneous subspaces $\A_n$ be invariant 
with respect to $R(\kappa)$ in such a way that $\forall \kappa \in \cK$ 
\begin{equation}
R(\kappa)\vert_{\A_0} := f((\ad \kappa)_0),
\quad
R(\kappa)\vert_{\A_n} := F((\ad \kappa)_n) \quad 
\forall n\in\ \cZ\setminus \{0\}.
\label{2.7}\end{equation}
For $n\in \cZ$ for which $\mathrm{dim}(\A_n)\neq 0$, 
these finite-dimensional linear operators are given similarly to (\ref{2.1}).
The assumption (\ref{2.6}) guarantees that the spectra $\sigma_\kappa^n$ 
do not intersect the poles of the corresponding 
meromorphic functions  $f$ and $F$ in (\ref{1.1}), whereby $R(\kappa)$ is well defined
for $\kappa\in \cK$.
If $\mathrm{dim}(\A_n)=0$, then $R(\kappa)\vert_{\A_n}$ is of course
understood to be the zero linear operator.  
Somewhat informally, we summarize (\ref{2.7}) 
by saying that $R(\kappa)$ equals $f(\ad \kappa)$ on $\K$ and
$F(\ad \kappa)$ on $\K^\perp$.

\bigskip\noindent 
{\bf Theorem 1.}  
{\em Let $\A$ be a graded, self-dual, complex Lie algebra satisfying the assumptions
given by (\ref{2.2})--(\ref{2.4}).
Take $\K:= \A_0$ and 
suppose the existence a nonempty, open domain  $\check \K\subset \K$
for which (\ref{2.6}) holds.
Then the $r$-matrix  $R: \check \K \rightarrow \mathrm{End}(\A)$ 
defined by (\ref{2.7}) satisfies the CDYBE (\ref{1.3}). 
Moreover, $R(\kappa)$ is an antisymmetric operator $\forall \kappa\in \check \K$,
and the $\K$-equivariance condition (\ref{1.6}) holds. }

\bigskip\noindent
{\bf Proof.}
Since the CDYBE (\ref{1.3}) is linear in $X, Y\in \A$, 
it is enough to verify it case by case for all possible choices of 
homogeneous elements $X$ and $Y$.
As a preparation, let us write the function $F$ in (\ref{1.1}) as
\begin{equation}
F(z)= \frac{1}{2} \frac{Q_+(z)}{ Q_-(z)} \quad
\hbox{with}\quad 
Q_\pm (z)= e^{\frac{z}{2}} \pm e^{-\frac{z}{2}},
\label{2.8}\end{equation}
and define the linear operators $Q_\pm(\kappa)$ on $\A$ by
\begin{equation}
Q_{\pm}(\kappa)= e^{K} \pm e^{-K}
\quad\hbox{with}\quad K:= \frac{1}{2} \ad \kappa
\quad\forall \kappa\in \check \K.
\label{2.9}\end{equation}
$Q_\pm(\kappa)$ are well defined operators on $\A$ since
their restrictions to any $\A_n$  are obviously well defined.
It follows from the definitions of the domain $\check \K$ 
and that of $R(\kappa)$ that $Q_-(\kappa)$ is an invertible
operator on $\A_n$ for any $n\neq 0$ and that we have
\begin{equation}
R(\kappa) Q_-(\kappa)= Q_-(\kappa) R(\kappa) =
\frac{1}{2} Q_+(\kappa)
\quad
\hbox{on} \quad\A_n \quad \forall n\neq 0.
\label{2.10}\end{equation}

We first consider the simplest case,  
\begin{equation}
X\in \A_m,\quad
 Y\in \A_n,
\quad
m\neq 0, \quad n\neq 0,\quad (m+n)\neq 0,
\label{2.11}\end{equation}
for which the derivative terms drop out from (\ref{1.3}).
Without loss of generality, we can now write 
\begin{equation}
X= Q_-(\kappa) \xi,
\qquad
Y= Q_-(\kappa) \eta 
\label{2.12}\end{equation}
with some $\xi\in \A_m$, $\eta \in \A_n$.
If we multiply (\ref{1.3}) from the left by the invertible operator
$4 Q_-(\kappa)$ on $\A_{m+n}$, then by using (\ref{2.10})
the required statement becomes
\begin{eqnarray}
&&
 Q_-(\kappa) [ Q_-(\kappa) \xi, Q_-(\kappa) \eta]+
Q_-(\kappa) [ Q_+(\kappa) \xi, Q_+(\kappa) \eta]
\nonumber\\
&&\qquad - Q_+(\kappa) \left( [ Q_-(\kappa)\xi, Q_+(\kappa) \eta]
+ [ Q_+(\kappa) \xi, Q_-(\kappa) \eta] \right)=0.
\label{2.13}\end{eqnarray}
We further spell out this equation by using 
that $e^{\pm K}$ are Lie algebra automorphism, and
thereby (\ref{2.13}) is verified in a straightforward manner.

Second, let us consider the case for which   
\begin{equation}
X\in \A_0,
\quad
Y\in \A_n, \quad n\neq 0.
\label{2.14}\end{equation}
Then the derivative term
$(\nabla_X R)(\kappa)(Y)$ appears in equation (\ref{1.3}).
To calculate this, we need the holomorphic complex function $h$ given by
\begin{equation}
z \mapsto h(z):= \frac{e^z -1}{z}.
\label{2.15}\end{equation}
We recall (e.g.~\cite{SW}, page 35) that for a curve
$t \mapsto A(t)$
of finite-dimensional linear operators one has the identity
\begin{equation}
\frac{ d e^{\pm A(t)}}{dt}= \pm e^{\pm A(t)} 
h(\mp \mathrm{ad}_{A(t)}) (\dot{A}(t)), 
\qquad
\dot{A}(t):= \frac{ dA(t)}{ dt}.
\label{2.16}\end{equation}
The right hand side of the above equation is defined by means  of 
the  Taylor expansion of $h$ around $0$, and of course   
\begin{equation}
(\mathrm{ad}_{A(t)})^j (\dot{A}(t))= [A(t), 
(\mathrm{ad}_{A(t)})^{j-1} (\dot{A}(t))],
\quad
j\in {\bN}, 
\quad
(\mathrm{ad}_{A(t)})^0 (\dot{A}(t))= \dot{A}(t).
\label{2.17}\end{equation}
In our case we consider the curve of linear operators on $\A_n$ given by 
\begin{equation}
t \mapsto \ad \kappa + t (\ad X).
\label{2.18}\end{equation}
Then (\ref{2.16}) leads to the formula
\begin{equation}
\left(\nabla_X e^{\pm K}\right)(Y)= \pm \frac{1}{2} e^{\pm K} [ h(\mp K) X, Y], 
\label{2.19}\end{equation}
where $K= \frac{1}{2}\ad \kappa$.  From this,  by taking the 
derivative of the identity 
$2 Q_- R = Q_+$ on $\A_n$ along the curve (\ref{2.18}) at $t=0$, 
we obtain 
\begin{equation}
4 Q_-(\kappa) (\nabla_X R)(\kappa) Y=
e^K [ h(-K)X, Y- 2 R(\kappa) Y] - e^{-K} [ h(K)X, Y+ 2 R(\kappa) Y].
\label{2.20}\end{equation} 
On the other hand, for (\ref{2.14})  the CDYBE (\ref{1.3}) is equivalent to
\begin{eqnarray}
&&4 Q_-(\kappa) (\nabla_X R)(\kappa) Y=
Q_-(\kappa) [ X, Y] + 4 Q_-(\kappa) [ R(\kappa) X, R(\kappa)Y] 
\nonumber\\
&&\qquad\qquad
- 2 Q_+(\kappa) \left( [X, R(\kappa) Y] + [ R(\kappa)X, Y]\right).
\label{2.21}\end{eqnarray}
We fix $\kappa\in \check \K$ arbitrarily,  
and write $Y= Q_-(\kappa) \eta $ with some 
 $\eta\in \A_n$. 
Then by a straightforward calculation, using that
$e^{\pm K}$ are Lie algebra automorphisms and collecting terms,
we obtain that the required equality of the right hand sides of
the last two equations is equivalent to
\begin{equation}
\left[\left( e^K h(-K) + e^{-K} h(K) - e^K  - e^{-K}\right) X, \eta\right]
= 2 \left[ \left( e^{-K} R(\kappa) - e^K R(\kappa)\right) X, \eta \right]. 
\label{2.22}\end{equation} 
Here $R(\kappa) X= f(2K)X$ with (\ref{1.2}), and the statement follows from the
equality of the corresponding complex analytic functions, namely   
\begin{equation}
e^z \frac{1-e^{-z}}{z} + e^{-z}\frac{e^z -1}{z} 
- e^z - e^{-z}=
 e^{-z} \left( \coth z - \frac{1}{z}\right) - e^z \left(\coth z - \frac{1}{z}\right),
\label{2.23}\end{equation}
which is checked in the obvious way. 

The third case to deal with is that of
\begin{equation}
X\in \A_{-n},
\quad Y\in \A_n, 
\quad n\neq 0,
\label{2.24}\end{equation}
for which the derivative term
$\langle X, (\nabla R)(\kappa) Y\rangle$ occurs in (\ref{1.3}).
At any fixed $\kappa \in \check \K$, 
we may write 
\begin{equation}
X= Q_-(\kappa) \xi,
\qquad
Y= Q_-(\kappa) \eta
\label{2.25}\end{equation}
with some $\xi\in \A_{-n}$, $\eta \in \A_n$.
We introduce the holomorphic function
\begin{equation}
z\mapsto g(z):= \frac{e^z - e^{-z}}{z},
\label{2.26}\end{equation}
and define $g(K)$ by the Taylor series of $g(z)$ around $z=0$.  
Then we can calculate that
\begin{equation}
\langle X, (\nabla R)(\kappa)Y\rangle = \frac{1}{2} g(K) [\eta, \xi].
\label{2.27}\end{equation}
To obtain this, note that 
\begin{equation}
\langle X, (\nabla R)(\kappa)Y\rangle = 
T^i \langle X, (\nabla_{T_i} R)(\kappa)Y\rangle 
\label{2.28}\end{equation}
with  dual bases $T_i$ and $T^i$ of $\A_0$, where 
$(\nabla_{T_i} R)(\kappa)Y$ is determined by (\ref{2.20}). 
By using these and  the invariance of the scalar product of $\A$, 
it is not difficult to rewrite (\ref{2.28}) in the form (\ref{2.27}).  
As for the non-derivative terms in  (\ref{1.3}), 
 with $X,Y$ in (\ref{2.25}) we find 
\begin{eqnarray}
&&[R(\kappa) X, R(\kappa) Y] - R(\kappa)\left(
[ X, R(\kappa) Y] + [R(\kappa) X, Y]\right) + \frac{1}{4} [X,Y]=
\nonumber\\
&&\qquad\qquad
\frac{1}{2} \left( Q_+(\kappa) - 2 R(\kappa)  Q_-(\kappa)\right) [\xi, \eta].
\label{2.29}\end{eqnarray}
It is easy to check that the sum of the right hand 
sides of (\ref{2.27}) and (\ref{2.29}) 
is zero, which finishes the verification of the CDYBE (\ref{1.3}) 
in the case (\ref{2.24}).

The remaining case is that of $X, Y\in \A_0$.
Then the variable $\kappa$ as well as 
all terms in (\ref{1.3}) lie in the subalgebra $\A_0$, 
and it is known \cite{EV,AM,BFP} that the formula
$\kappa \mapsto f(\ad \kappa)$ (\ref{1.2}) 
defines a solution of the CDYBE on any finite-dimensional 
self-dual Lie algebra. This completes the verification of
the CDYBE (\ref{1.3}).

The antisymmetry of $R(\kappa)$ follows from
(\ref{2.7}) since $\ad \kappa$ is antisymmetric by the invariance of 
$\langle\ ,\ \rangle$ and both $f$ and $F$ are odd functions. 
Finally, the equivariance property (\ref{1.6}) is also 
easily verified from (\ref{2.7}) by using that 
for any finite-dimensional 
linear operator given by (\ref{2.1}) one has 
\begin{equation}
{\frac{d H(A(t))}{dt}}\vert_{t=0} =
 \frac{1}{2\pi i} \oint_\Gamma dz H(z) (z I_V - A)^{-1} 
\dot{A}(0) (z I_V - A)^{-1}
\label{2.30}\end{equation}
along any smooth curve $t\mapsto A(t)$ for which $A(0)=A$.
{\em Q.E.D.}

\bigskip

We conclude this section by describing the
tensorial interpretation of the CDYBE (\ref{1.3}) for 
the $r$-matrices of theorem 1.
For this,  consider dual bases  
$T_i[n]$ and $T^j[n]$ of $\A$ 
($n\in \cZ, i,j=1,\ldots,\mathrm{dim}(\A_n)$),
which satisfy $T_i[n]\in \A_n$ and 
$ \langle T_i[m], T^j[n]\rangle = \delta_{m, -n} \delta_i^j$.
Then introduce  $r^\pm :\check \K \rightarrow \A \otimes \A$ 
by\footnote{Here $\A\otimes \A$ is a 
completion of the algebraic tensor product 
containing the elements that are associated with the 
linear operators on $\A$.}  
\begin{equation}
r^{\pm}(\kappa) := \sum_{n\in\cZ} \sum_{i=1}^{\mathrm{dim}(\A_n)}
\left( (R(\kappa) T_i[n]) \otimes T^i[-n] \pm \frac{1}{2} 
T_i[n] \otimes T^i[-n] \right).
\label{2.31}\end{equation}
In fact, 
as a consequence of the properties of $R$ established in theorem 1,
$r^\pm$ satisfies the tensorial version of the CDYBE given by 
\begin{eqnarray}
&& [ r^{s}_{12}(\kappa), r^{s}_{13}(\kappa)] 
+[ r^{s}_{12}(\kappa), r^{s}_{23}(\kappa)] +
[ r^{s}_{13}(\kappa), r^{s}_{23}(\kappa)]\qquad \nonumber\\
&& \quad  + 
T_{j}[0]^1 \frac{\partial}{\partial \kappa_j} r_{23}^{s}(\kappa) 
-T_{j}[0]^2 \frac{\partial}{\partial \kappa_j}  r_{13}^{s}(\kappa) 
+T_{j}[0]^3 \frac{\partial}{\partial \kappa_j}  r_{12}^{s}(\kappa)  =0,
\quad s=\pm,\quad 
\label{2.32}\end{eqnarray}
where $\kappa_j:= \langle \kappa,T_j[0]\rangle$.
Here the standard notations are used, 
$T_j[0]^1:= T_j[0] \otimes 1\otimes 1$, 
$r_{12}^s := r^s \otimes 1$  etc.
The expression on the 
left hand side of (\ref{2.32}) belongs to 
a completion of $\A\otimes \A\otimes \A$; 
it has a unique expansion in the basis 
$T_{i_1}[n_1] \otimes T_{i_2}[n_2]\otimes T_{i_3}[n_3]$
of $\A\otimes \A\otimes \A$. 
Similarly to the CYBE, the CDYBE
(\ref{2.32}) is  compatible with homomorphisms of $\A$. 
This means that if $\pi_i : \A \rightarrow \G^i$ ($i=1,2,3$)
are (possibly different) 
homomorphisms of $\A$ into (possibly different) Lie algebras $\G^i$,
then we can obtain a $\G^1 \otimes \G^2 \otimes \G^3$-valued
equation from (\ref{2.32}) by the obvious application of the map 
$\pi_1 \otimes \pi_2 \otimes \pi_3$ to all objects on the
left hand side of (\ref{2.32}).  
More precisely, to take into account the unit element $1$,
here one uses the extensions of these Lie algebra homomorphisms to the
corresponding universal enveloping algebras.  

\section{Applications to affine Lie algebras}
\setcounter{equation}{0}

Let  $\G$ be a finite-dimensional complex, self-dual Lie algebra equipped 
with an invariant `scalar product' denoted as $B: \G \times \G \rightarrow \bC$.
Suppose that $\mu$ is a finite order automorphism of $\G$ 
that preserves the bilinear form $B$ and has nonzero fixed 
points\footnote{The last two properties are automatic if $\G$ 
is simple or $\mu=\mathrm{id}$, which are included as special cases.}.
With this data, one may associate the twisted loop algebra $\ell(\G,\mu)$ 
and the affine Lie algebra $\A(\G,\mu)$ obtained by adding the natural 
derivation to the central extension of $\ell(\G,\mu)$.
We below show that theorem 1 is directly applicable to $\A(\G,\mu)$. 
Then we explain that the resulting dynamical $r$-matrices on 
$\A(\G,\mu)$ admit a reinterpretation as one-parameter
families of $r$-matrices on $\ell(\G,\mu)$.
By applying evaluation homomorphisms to the corresponding 
$\ell(\G,\mu)\otimes \ell(\G,\mu)$-valued $r$-matrices, 
we finally derive spectral-parameter-dependent $\G\otimes \G$-valued
dynamical $r$-matrices.
These results were announced in \cite{FGP1} without presenting proofs, 
which are provided here.

\subsection{Application of theorem 1 to  $\A(\G,\mu)$}

Any automorphism $\mu$  of order $N$, $\mu^N= \mathrm{id}$,
gives rise to a decomposition of $\G$ as   
\begin{equation}
\G= \oplus_{a\in \E_\mu} \G_a,
\qquad
\E_\mu\subset \{ 0, 1, \ldots, (N-1)\,\},
\label{3.1}\end{equation}
with the eigensubspaces 
\begin{equation}
\G_a:= \{ \xi\in \G\,\vert\, \mu(\xi)= 
\exp(\frac{ia2\pi}{N}) \xi\,\} \neq \{ 0\}.
\label{3.2}\end{equation}
Since we assumed that $B(\mu \xi, \mu \eta)=B(\xi,\eta)$ ($\forall \xi,\eta\in \G$),
$\G_a$ is perpendicular to $\G_b$ with respect to the form $B$
unless $a+b=N$ or $a=b=0$.
This implies that if a nonzero $a$
belongs to the index set $\E_\mu$ then so does $(N-a)$.
We assume that  $0 \in \E_\mu$, and thus
$\G_0\neq \{ 0\}$ is a self-dual subalgebra of $\G$.

The twisted (or untwisted if we choose $\mu=\mathrm{id}$) 
loop algebra $\ell(\G,\mu)$ is 
the subalgebra of $\G\otimes {\bC}[t,t^{-1}]$
generated by the elements of the form
\begin{equation}
\xi^{n_a}:= \xi \otimes t^{n_a}
\quad\hbox{with}\quad
\xi\in \G_a,
\quad
n_a = a + m N,
\,\,\, m\in {\bZ},
\label{3.3}\end{equation}
where $t$ is a formal variable.
The `affine Lie algebra' $\A(\G,\mu)$ is then  
introduced as 
\begin{equation}
\A(\G,\mu):= \ell(\G, \mu) 
\oplus {\bC} d \oplus {\bC} \hat c 
\label{3.4}\end{equation}
with the Lie bracket of its generators defined by  
\begin{equation}
[\xi^{n_a}, \eta^{p_b}] = [\xi,\eta]^{n_a + p_b}
+ n_a \delta_{n_a,- p_b} B(\xi,\eta)\hat c,
\quad
\forall \xi\in \G_a,\,\,\,\eta\in \G_b,
\label{3.5}\end{equation}
\begin{equation}
[ d, \xi^{n_a}]= n_a \xi^{n_a},
\quad
[\hat c, d] = [\hat c, \xi^{n_a}]=0. 
\label{3.6}\end{equation}
$\A(\G,\mu)$ is a self-dual Lie algebra as it carries the 
scalar product $\langle\ ,\ \rangle$ given by
\begin{equation}
\langle \xi^{n_a}, \eta^{p_b}\rangle  = 
\delta_{n_a, -p_b} B(\xi,\eta),
\quad
\langle \hat c, d \rangle =1,
\quad
\langle d, \xi^{n_a}\rangle = \langle \hat c, \xi^{n_a}\rangle =0.
\label{3.7}\end{equation}
We obtain a $\bZ$-gradation of $\A(\G,\mu)$ by the decomposition  
\begin{equation} 
\A(\G,\mu)= \oplus_{n\in (\E_\mu + N {\bZ})} \A(\G,\mu)_n=\oplus_{n\in \bZ} \A(\G,\mu)_n,
\label{3.8}\end{equation}
where $\A(\G,\mu)_n$ is the eigensubspace of $\ad d$ with eigenvalue $n$ if
$n\in (\E_\mu + N {\bZ})$, and $\A(\G,\mu)_n=\{0\}$ 
if $n\notin (\E_\mu + N {\bZ})$. 
We need to introduce these zero subspaces for notational consistency, 
since  $(\E_\mu + N {\bZ})$ is not necessarily a group in general.
This is also consistent with the fact that (\ref{3.5}) 
gives zero if $(n_a+ p_b)\notin (\E_\mu + N {\bZ})$.
The gradation given by (\ref{3.8}) clearly   
satisfies equations (\ref{2.2})--(\ref{2.4}), 
where now $\cZ :=\bZ$.
We below regard $\G_0$ as a subspace of $\A(\G,\mu)$ by
identifying $\xi\in \G_0$ with $\xi\otimes t^0\in \A(\G,\mu)$,
whereby we can write 
\begin{equation}
 \A(\G,\mu)_0 = \G_0 \oplus {\bC} d \oplus {\bC}\hat c.
\label{3.9}\end{equation}

Since we wish to apply theorem 1, we now set $\A:=\A(\G,\mu)$ and 
$\K:= \A(\G,\mu)_0$.
We parametrize the general element $\kappa\in \K$ as
\begin{equation}
\kappa= \omega + k d + l \hat c,
\qquad
\omega\in \G_0,\quad k,l\in {\bC}.
\label{3.10}\end{equation}
It follows from the above that 
formula (\ref{2.7}) provides us with a dynamical $r$-matrix 
$R: \check \K \rightarrow \mathrm{End}(\A)$ if we can find 
a nonempty, open domain $\check \K\subset \K$ whose elements satisfy 
the conditions given in (\ref{2.6}).
The point is that  we can indeed find such a domain,
and actually the maximal domain has the form 
\begin{equation}
\check \K= \{ \,\kappa = \omega + k d + l \hat c\,\vert\,
l\in {\bC}, \,\,
k \in ({\bC}\setminus {\bR}i),
\,\,\,  \omega \in \B_k \},
\label{3.11}\end{equation}
where $\B_k\subset \G_0$ is described as follows.
Let $\lambda_a$ denote an eigenvalue of the operator 
$\ad \omega\vert_{\G_a}$ associated with $\omega\in \G_0$.
By definition, the subset $\B_k\subset \G_0$ consists of those $\omega\in \G_0$ 
whose eigenvalues satisfy  the following conditions:
\begin{equation}
(\lambda_a + k (a+m N))\notin 2\pi i \bZ
\quad
\forall m\in \bZ, \quad 
\forall a\in \E_\mu\setminus \{ 0\},
\label{3.12}\end{equation} 
\begin{equation}
\lambda_0\notin 2\pi i\bZ^*
\quad\hbox{and}\quad
(\lambda_0 + k m N)\notin 2\pi i \bZ
\quad
\forall m\in \bZ^*.
\label{3.13}\end{equation} 
If we note that for $\xi^{n_a}$ in (\ref{3.3}) and 
$\kappa \in \K$ written as in (\ref{3.10}) one has 
\begin{equation}
(\ad \kappa)( \xi^{n_a}) = k n_a \xi^{n_a} + [\omega, \xi]^{n_a},  
\label{3.14}\end{equation}
then 
the conditions in (\ref{3.12}) and (\ref{3.13}) are recognized to be the translation 
of the condition in (\ref{2.6}) to our case.
The set $\check \K$ defined by these requirements obviously contains 
the elements of the form $\kappa = kd + l\hat c$ for any $k\in (\bC\setminus i \bR)$,
$l\in \bC$, and therefore it is nonempty.
It is not difficult to see that $\check \K \subset \K$ in (\ref{3.11}) is an open subset,
for which one needs $k$ to have a nonzero real part, and 
$\B_k \subset \G_0$ is a nonempty open subset as well.
For completeness, we present a proof of these statements 
in appendix A.
   
\subsection{One-parameter family of $r$-matrices on $\ell(\G,\mu)$}

We now reinterpret the dynamical $r$-matrices
$R: \check \K \rightarrow \mathrm{End}(\A(\G,\mu))$ 
constructed in subsection 3.1 as a family of $r$-matrices
\begin{equation}
R_k:  \B_k \rightarrow \mathrm{End}(\ell(\G,\mu)),
\label{3.15}\end{equation}
where the parameter $k$ varies in 
$(\bC\setminus i \bR)$
and the $k$-dependent domain $\B_k\subset \G_0$ appears in (\ref{3.11}).
For any $\omega \in \B_k$, the operator $R_k(\omega)$ is given by
\begin{equation}
R_k(\omega)\eta:= f(\ad \omega)\eta, 
\qquad
R_k(\omega) \xi^{n_a}:=
F( k n_a +  \ad \omega) \xi^{n_a} 
\label{3.16}\end{equation} 
$\forall \eta \in \G_0 =\ell(\G,\mu)_0$ and 
$\forall \xi^{n_a}\in\ell(\G,\mu)_{n_a}$ with $n_a\neq 0$.
In other words, by regarding $\ell(\G,\mu)$ as a subspace
of $\A(\G,\mu)$, we have  $R_k(\omega) X= R(\kappa) X$
for $X\in \ell(\G,\mu)$ and $\kappa\in \check \K$.

It is an easy consequence of theorem 1 that
$R_k$ satisfies the operator version  of the CDYBE for any fixed $k$:
\begin{eqnarray}
&&[ R_k X, R_k Y] -R_k( [X, R_k Y]+ [R_k X,Y]) 
+\langle X, (\nabla R_k) Y\rangle  \nonumber\\
&& \quad + (\nabla_{Y_0} R_k) X -
(\nabla_{X_0} R_k) Y =-\frac{1}{4} [X,Y], 
\qquad \forall X, Y\in \ell(\G,\mu).
\label{3.17}\end{eqnarray}
Here the Lie brackets are evaluated in $\ell(\G,\mu)$, 
$X_0$ is the grade $0$ part of $X$, and the scalar product 
on $\ell(\G,\mu)$ is given by the restriction of (\ref{3.7}).
This equation is verified by a simplified version of the 
calculation done in the proof of theorem 1,
the simplification being that $\hat c$ has now been set to zero.
It is also clear that 
$R_k: \B_k \rightarrow \mathrm{End}(\ell(\G,\mu))$
is a $\G_0$-equivariant map in the natural sense.

For later purpose, we here introduce the shifted $r$-matrices
\begin{equation}
R_k^{\pm}:= R_k \pm \frac{1}{2} I,
\label{3.18}\end{equation}
where $I$ is the identity operator on $\ell(\G,\mu)$.
By using the scalar product, we associate with these operator 
valued maps the corresponding 
$\ell(\G,\mu)\otimes \ell(\G,\mu)$-valued maps.
These are denoted respectively as 
\begin{equation}
r^{k,\pm}: \B_k \rightarrow \ell(\G,\mu) \otimes \ell(\G,\mu).
\label{3.19}\end{equation} 
By translating the CDYBE into tensorial terms, 
(\ref{3.17}) becomes  
\begin{eqnarray}
&& [ r^{k,s}_{12}(\omega), r^{k,s}_{13}(\omega)] 
+[ r^{k,s}_{12}(\omega), r^{k,s}_{23}(\omega)] +
[ r^{k,s}_{13}(\omega), r^{k,s}_{23}(\omega)] \nonumber\\
&& \qquad  + 
T_{j}^1 \frac{\partial}{\partial \omega_j} r_{23}^{k,s}(\omega) 
-T_{j}^2 \frac{\partial}{\partial \omega_j}  r_{13}^{k,s}(\omega) 
+T_{j}^3 \frac{\partial}{\partial \omega_j}  r_{12}^{k,s}(\omega)  =0,
\quad s=\pm,\quad 
\label{3.20}\end{eqnarray}
where $\omega_j:= B( \omega,T_j)$ with a basis $T_j$ of $\G_0$.

\subsection{Spectral-parameter-dependent $r$-matrices}

The loop algebra $\ell(\G,\mu)$ admits an `evaluation homomorphism' 
$\pi_v:\ell(\G,\mu)\rightarrow \G$ for any fixed $v\in {\bC}^*$,
\begin{equation}
\pi_v: \xi \otimes t^n \mapsto v^n \xi  
\qquad \forall (\xi\otimes t^n) \in \ell(\G,\mu).
\label{3.21}\end{equation}
It is well known that spectral-parameter-dependent 
$\G\otimes \G$-valued $r$-matrices 
 may be obtained by applying these homomorphisms to 
$\ell(\G,\mu)\otimes \ell(\G,\mu)$-valued $r$-matrices.
In the context of dynamical $r$-matrices, Etingof and Varchenko \cite{EV}
used this method to derive Felder's elliptic dynamical $r$-matrices from 
the basic trigonometric dynamical $r$-matrices of the 
(untwisted) affine Kac-Moody Lie algebras. 
We here apply the same procedure to the general family  
of dynamical $r$-matrices introduced in eq.~(\ref{3.19}).
As for the presentation below, we find it convenient 
to first provide a self-contained definition 
of the spectral-parameter-dependent $r$-matrices and 
show afterwards how they are obtained from the evaluation 
homomorphisms.   

We start by collecting some meromorphic functions 
and identities that will be useful. 
Consider the standard theta function\footnote{We have 
$\theta_1(z\vert \tau)=\vartheta_1(\pi z\vert \tau)$ 
with $\vartheta_1$ in \cite{WW}.} 
\begin{equation}
\theta_1(z\vert\tau):= -\sum_{j\in {\bZ}} 
\exp\left( \pi i (j+ \frac{1}{2})^2 \tau + 
2\pi i (j+\frac{1}{2}) (z + \frac{1}{2})\right),
\quad 
\Im(\tau) >0,
\label{3.22}\end{equation}
which is holomorphic on $\bC$ and has simple zeros at the points 
of the lattice
\begin{equation}
\Omega:= \bZ + \tau \bZ.
\end{equation}
Recall that $\theta_1$ is odd in $z$ and satisfies 
\begin{equation}
\theta_1(z+1\vert \tau)=-\theta_1(z\vert \tau),
\quad
\theta_1(z+\tau\vert \tau)=-q^{-1}e^{-2\pi i z} \theta_1(z\vert \tau),
\quad
q:= e^{\pi i \tau}. 
\label{3.23}\end{equation}
Define now the function
\begin{equation}
\chi(w,z\vert \tau):= \frac{1}{2\pi i} 
\frac{\theta_1(\frac{w}{2\pi i}+z\vert\tau)\theta'_1(0\vert\tau)}{\theta_1(z\vert \tau) 
\theta_1(\frac{w}{2\pi i}\vert\tau)}.
\label{3.24}\end{equation}
This function is holomorphic in $w$ and in $z$ at the points 
\begin{equation}
(w,z)\in (\bC\setminus 2\pi i \Omega)\times (\bC \setminus \Omega). 
\label{3.25}\end{equation}
The following important identity holds:
\begin{equation}
\chi(w,z) = \frac{1}{2}
\sum_{n\in {\bZ}} e^{2\pi i z n} 
\left[ 1 + \coth(\frac{w}{2} + \pi i \tau n)\right] 
\label{3.26}\end{equation}
on the domain 
\begin{equation}
D:= \{(w,z)\,\vert\, w\in (\bC\setminus 2\pi i \Omega),\,\,\,
-\Im(\tau) <\Im(z) <0\,\}.
\label{3.27}\end{equation}
All terms in the sum are holomorphic on $D$, 
the convergence is absolute at any fixed $(w,z)\in D$, and is uniform
on compact subsets of $D$.
The verification of (\ref{3.26}) is a routine matter,
example 13 on page 489 of \cite{WW} contains a closely related statement.

We also need the functions 
\begin{equation}
\chi_a(w,z\vert \tau):= 
e^{\frac{2\pi i a z}{N}}\left(\chi(w+2\pi i\frac{ a}{N}\tau, z\vert\tau) 
-\frac{\delta_{a,0}}{w}\right),
\label{3.28}\end{equation}
where $a\in \{0,1,\ldots, (N-1)\}$ with some positive integer $N$.
The function $\chi_{a}(w,z\vert\tau)$ is 
holomorphic in $w$ and in $z$ if $(w,z)$ belongs to the domain 
\begin{equation}
(\bC\setminus 2\pi i \Omega_a)\times (\bC \setminus \Omega)
\quad
\hbox{where}\quad
\Omega_a:= \left(\Omega - \frac{a}{N}\tau\right) \setminus \{ 0\}.
\label{3.29}\end{equation}
By using the notation
\begin{equation}
f_a(w):=
 \frac{1}{2}\left[1+\coth\frac{w}{2}\right] - \frac{\delta_{a,0}}{w},
\label{3.30}\end{equation}
we have the identity
\begin{equation}
\chi_a(w,z\vert\tau) =e^{\frac{2\pi i a z}{N}}\left(f_a(w+ 2\pi i \frac{a}{N}\tau)
+\frac{1}{2} \sum_{n\in {\bZ}^*} e^{2\pi i z n} 
\left[ 1 + \coth(\frac{w}{2} +\pi i \frac{a}{N}\tau+ \pi i \tau n)\right] 
\right)
\label{3.31}\end{equation}
on the domain 
\begin{equation}
D_a:= \{(w,z)\,\vert\, w\in (\bC\setminus 2\pi i \Omega_a),\,\,\,
- \Im(\tau) <\Im(z) <0\,\}
\label{3.32}\end{equation}
for any $a\in \{ 0,1,\ldots, (N-1)\}$.
All terms in the sum are holomorphic on $D_a$, 
the convergence is absolute at any $(w,z)\in D_a$, and is uniform
on compact subsets of $D_a$. 

Let now $\mu$ be an automorphism of $\G$ of order $N$
as considered previously and fix $\tau$ with $\Im(\tau)>0$.
For any $\omega\in \G_0$ and $a\in \E_\mu$, 
let $\sigma((\ad\omega)_a)$ be the spectrum of the linear operator
$(\ad\omega)_a:=\ad \omega\vert_{\G_a}$.
Define $\B^\tau\subset \G_0$ by 
\begin{equation}
\B^\tau:= \{\, \omega \in \G_0\,\vert\,
\sigma((\ad\omega)_a)\cap 2\pi i \Omega_a =\emptyset\quad
\forall a\in \E_\mu\,\}.
\label{3.33}\end{equation}
It is easy to verify that 
\begin{equation}
\B^\tau = \B_k \quad\hbox{if}\quad
\tau= \frac{kN}{2\pi i},
\label{3.34}\end{equation}
where $\B_k\subset \G_0$ appears in (\ref{3.11}).
In particular, $\B^\tau$ is an open subset of $\G_0$ 
that contains the origin. 
By using the above notations, we now  define the function $\R_\tau$ as
\begin{equation}
\R_\tau: \B^\tau \times (\bC\setminus \Omega) \rightarrow \mathrm{End}(\G),
\qquad
\R_\tau(\omega,z)\vert_{\G_a}:= \chi_a((\ad \omega)_a, z\vert \tau).
\label{3.35}\end{equation} 
It follows from the properties of the holomorphic functional
calculus on Banach algebras \cite{DS} that $\R_\tau$ is
well defined and is holomorphic in its variables.
Next we introduce also the holomorphic function
\begin{equation}
r^\tau: \B^\tau \times (\bC\setminus \Omega) \rightarrow 
\G\otimes \G,
\qquad
r^\tau(\omega,z):= B(T_\alpha, \R_\tau(\omega,z)T_\beta) 
T^\alpha \otimes T^\beta,
\label{3.36}\end{equation}     
where $T_\alpha$, $T^\beta$ are dual bases of $\G$.
We now state one of our main results. 

\bigskip\noindent
{\bf Proposition 2.}
{\em 
The function $r^\tau$ introduced above satisfies the
spectral-parameter-dependent version of the CDYBE:
\begin{eqnarray}
&& [ r^{\tau}_{12}(\omega, z_{12}), r^{\tau}_{13}(\omega, z_{13})] 
+[ r^{\tau}_{12}(\omega,z_{12}), r^{\tau}_{23}(\omega,z_{23})] +
[ r^{\tau}_{13}(\omega,z_{13}), r^{\tau}_{23}(\omega,z_{23})] \nonumber\\
&& \qquad  + 
T_{j}^1 \frac{\partial}{\partial \omega_j} r_{23}^{\tau}(\omega,z_{23}) 
-T_{j}^2 \frac{\partial}{\partial \omega_j}  r_{13}^{\tau}(\omega,z_{13}) 
+T_{j}^3 \frac{\partial}{\partial \omega_j}  r_{12}^{\tau}(\omega,z_{12})  
=0,
\label{3.37}\end{eqnarray} 
where $z_{\alpha\beta}= (z_\alpha -z_\beta)\in (\bC\setminus \Omega)$,
$\omega\in \B^\tau$, and 
$\omega_j:= B( \omega,T_j)$ with a basis $T_j$ of $\G_0$.
Furthermore, $r^\tau$ has the properties
\begin{equation}
\mathrm{Res}_{z=0} r^\tau(\omega,z)= \frac{1}{2\pi i} T^\alpha \otimes T_\alpha,
\qquad\qquad 
\left(r^\tau(\omega,z)\right)^T + r^\tau(\omega,-z)=0,
\label{3.38}\end{equation}
where $\left(r^\tau(\omega,z)\right)^T:=  
B(T_\alpha, \R_\tau(\omega,z)T_\beta) T^\beta \otimes T^\alpha$
with dual bases $T_\alpha$, $T^\beta$ of $\G$, and 
\begin{equation}
\frac{d}{dx} 
r^\tau(e^{\ad Tx}(\omega),z)\vert_{x=0}= 
[T\otimes 1 + 1 \otimes T, r^\tau(\omega,z)]
\qquad
\forall T\in \G_0.
\label{extra1}\end{equation}
}
\bigskip

The statements in (\ref{3.38}) follow immediately  from 
the definition (\ref{3.35}), (\ref{3.36}) and 
the properties of the meromorphic functions $\chi_a$ in (\ref{3.28}).
For the first equality in (\ref{3.38}), one can check that 
\begin{equation}
\mathrm{Res}_{z=0} \chi_a(w,z\vert \tau)= \frac{1}{2\pi i},
\qquad 0\leq a<N.
\end{equation}
For the second statement, one uses the invariance 
of the scalar product $B$ on $\G$ and 
\begin{equation} 
\chi_0(-w,z\vert\tau)=-\chi_0(w,-z\vert\tau),
\quad 
\chi_{a}(-w,z\vert\tau)=-\chi_{N-a}(w,-z\vert\tau),
\quad 
0<a <N.
\end{equation}
The $\G_0$-equivariance property (\ref{extra1}) is obvious from 
the definition of $r^\tau$.
As for the CDYBE (\ref{3.37}), it is consequence of the
following result.

\bigskip\noindent
{\bf Proposition 3.}
{\em The dynamical $r$-matrix $r^\tau$ 
given by (\ref{3.35}), (\ref{3.36})
results by evaluation homomorphism from the dynamical $r$-matrix
$r^{k,+}$ in (\ref{3.19}). 
More precisely, if we set   
\begin{equation}
\tau= \frac{kN}{2\pi i}\quad\hbox{and}\quad \frac{v_1}{v_2}= \exp(\frac{2\pi i z}{N}) 
\quad\hbox{with}\quad
- \Im(\tau) <\Im(z) <0,
\label{3.41}\end{equation}
then the evaluation homomorphism (\ref{3.21}) yields the
relation
\begin{equation}
(\pi_{v_1} \otimes \pi_{v_2})( r^{k,+}(\omega))= r^\tau(\omega,z)
\qquad
\forall \omega\in \B_k=\B^\tau.
\label{3.42}\end{equation}
}
\bigskip

\noindent{\bf Proof.}
The left hand side of (\ref{3.42}) gives only a formal
infinite sum in general. Below we first calculate this sum, 
and then notice that it converges to the function
on the right hand side of (\ref{3.42}) if the variables
satisfy (\ref{3.41}).

Let $T_{a, j}$ and $T_a^{j}$ ($j=1,\ldots, \mathrm{dim}(\G_a)$)
denote bases of $\G_a$ ($a\in \E_\mu$) subject 
to the relations
\begin{equation}
\langle T_{0,j}, T_0^l\rangle = \delta^l_j,
\qquad
\langle T_{a,j}, T_{N-a}^l\rangle = \delta^l_j,
\quad
\forall a\in \E_\mu\setminus \{0\}.
\label{3.43}\end{equation} 
Introduce corresponding bases of $\ell(\G,\mu)$:
\begin{equation}
T_{a,j}[n_a]:= T_{a,j}\otimes t^{n_a},
\qquad
T_a^j[n_a]:= T_a^j\otimes t^{n_a},
\quad
\forall a\in \E_\mu,\,\,n_a\in (a + N {\bZ}).
\label{3.44}\end{equation}
By definition, we then have
\begin{eqnarray}
&&r^{k,+}(\omega)= \sum_{j,l=1}^{\mathrm{dim}(\G_0)} 
\sum_{n_0\in N {\bZ}}
\langle T_{0,j}[-n_0], R_k^+(\omega) T_{0,l}[n_0]\rangle\,
T_0^j[n_0] \otimes T_0^l[-n_0]\nonumber\\
&&\,\,
+\sum_{a\in \E_\mu\setminus \{0\}}\sum_{j,l=1}^{\mathrm{dim}(\G_a)} 
\sum_{n_a\in (a+N {\bZ})}
\langle T_{N-a,j}[-n_a], R_k^+(\omega) T_{a,l}[n_a]\rangle\,
T_a^j[n_a] \otimes T_{N-a}^l[-n_a].\qquad\quad
\label{3.45}\end{eqnarray}
By substituting the definition of 
$R_k^+(\omega)$, (\ref{3.18}) with (\ref{3.16}),
we obtain that
\begin{equation}
\langle T_{N-a,j}[-n_a], 
R_k^+(\omega) T_{a,l}[n_a]\rangle = 
B( T_{N-a,j}, 
(   F(k n_a + \ad \omega) +\frac{1}{2}) T_{a,l} )
\label{3.46}\end{equation}
for $a\in \E_\mu\setminus \{0\}$, and 
\begin{equation}
\langle T_{0,j}[-n_0], 
R_k^+(\omega) T_{0,l}[n_0]\rangle = 
B( T_{0,j},  (F(k n_0 + \ad \omega) +\frac{1}{2}) T_{0,l} ),
\quad n_0\neq 0,
\label{3.47}\end{equation}
\begin{equation}
\langle T_{0,j}[0], 
R_k^+(\omega) T_{0,l}[0]\rangle = 
B( T_{0,j}, 
(f( \ad \omega) +\frac{1}{2}) T_{0,l} ),
\label{3.48}\end{equation}
where the functions $f$ and $F$ are given in (\ref{1.1}).
This implies that the left hand side of (\ref{3.42}) can  
be written in the following form: 
\begin{eqnarray}
&& (\pi_{v_1} \otimes \pi_{v_2})( r^{k,+}(\omega)) = 
\sum_{j,l=1}^{\mathrm{dim}(\G_0)} 
B( T_{0,j}, \psi_0((\ad \omega)_0, z\vert k) T_{0,l})\,
T_0^j \otimes T_0^l\nonumber\\
&&\quad
+\sum_{a\in \E_\mu\setminus \{0\}}\sum_{j,l=1}^{\mathrm{dim}(\G_a)} 
B( T_{N-a,j}, \psi_a((\ad \omega)_a, z\vert k) T_{a,l})\,
T_a^j \otimes T_{N-a}^l
\label{3.49}\end{eqnarray}
with
\begin{equation}
\psi_a((\ad \omega)_a, z\vert k) =
\frac{\exp(\frac{2\pi i a z}{N})}{2}
\sum_{m\in {\bZ}} e^{2\pi i z m} 
\left[1 + 
\coth \frac{k N m + ka + (\ad \omega)_a}{2}\right],
\quad a\neq 0,  
\label{3.50}\end{equation}
and 
\begin{equation}
\psi_0((\ad \omega)_0, z\vert k) =
\left[\frac{1}{2} + f((\ad \omega)_0)\right] +
\frac{1}{2}
\sum_{m\in {\bZ}^*} e^{2\pi i z m} 
\left[1 + 
\coth \frac{k N m  + (\ad \omega)_0}{2}\right].
\label{3.51}\end{equation}

To obtain the $a\neq 0$ terms in (\ref{3.49}) from (\ref{3.45}), 
we used (\ref{3.46}) and  the 
parametrization $\frac{v_1}{v_2}=\exp(\frac{2\pi i z}{N})$,
whereby  
\begin{eqnarray}
&&
\sum_{n_a\in (a+N {\bZ})}
\langle T_{N-a,j}[-n_a], R_k^+(\omega) T_{a,l}[n_a]\rangle\,
(\pi_{v_1}\otimes \pi_{v_2})\left( T_a^j[n_a] \otimes T_{N-a}^l[-n_a]\right)
\qquad\qquad
\nonumber\\ 
&& =  
\frac{1}{2} e^\frac{2\pi i a z}{N}
\sum_{m\in {\bZ} }e^{2\pi i z m} 
B(T_{N-a, j}, [ 1+ 2F(k a+kmN + \ad \omega)] T_{a,l}  )   T_a^j \otimes T_{N-a}^l
\nonumber\\
&& = 
B(T_{N-a, j}, 
\frac{1}{2} e^\frac{2\pi i a z}{N}
\sum_{m\in {\bZ} }e^{2\pi i z m} [1+2F(k a +kmN+ \ad \omega)] T_{a,l}  )   
T_a^j \otimes T_{N-a}^l.
\label{3.52}\end{eqnarray}
This leads to (\ref{3.49}) with (\ref{3.50}) by inserting the definition 
of $F$ (\ref{1.1}) and noting that $(\ad \omega) T_{a,l} =(\ad \omega)_a T_{a,l}$. 
The $a=0$ term is dealt with in a similar manner. 

Now we come to the main point. 
We notice that if on the right hand sides of (\ref{3.50}) and (\ref{3.51}) 
$(\ad \omega)_a$ is replaced by a complex variable $w$ and one uses also 
$\tau = \frac{kN}{2\pi i}$, then these series become precisely identical 
with the corresponding series in (\ref{3.31}), which are convergent on the 
domain $D_a$ (\ref{3.32}) for any $a\in \E_\mu$.
Since these are absolute convergent series and the convergence is uniform 
on compact subsets of $D_a$, it follows that the corresponding operator series
in (\ref{3.50}), (\ref{3.51}) converge, too.
Therefore, if 
\begin{equation}
\tau = \frac{kN}{2\pi i}, \qquad
\omega \in \B^\tau,\qquad 
 -\Im(\tau) < \Im(z) <0,
\label{3.53}\end{equation}
then $\psi_a( (\ad \omega)_a, z\vert k)\in \mathrm{End}(\G_a)$ is well defined 
by the corresponding series  in (\ref{3.50}), (\ref{3.51}),
and on this domain we obtain 
\begin{equation}
\psi_a((\ad\omega)_a, z\vert k)= 
\chi_a((\ad \omega)_a, z\vert \tau),
\qquad
\forall a\in \E_\mu.
\label{3.54}\end{equation}
If we now compare (\ref{3.49})  with the definition of $r^\tau$ given by 
(\ref{3.35}), (\ref{3.36}), 
then (\ref{3.54}) allows us to conclude 
that $(\pi_{v_1} \otimes \pi_{v_2})( r^{k,+}(\omega))= r^\tau(\omega,z)$
holds indeed on the domain given by (\ref{3.41}). 
{\em Q.E.D.} 
\bigskip

It is clear from the proof that (\ref{3.41}) is necessary for 
(\ref{3.42}); the series appearing in (\ref{3.50}) and (\ref{3.51}) 
do not converge  if $z$ lies outside the strip in (\ref{3.41}). 
Thus, by applying $\pi_{v_1}\otimes \pi_{v_2} \otimes \pi_{v_3}$ to 
the CDYBE (\ref{3.20}), proposition 3 directly implies proposition 2 
if $z_{12}$, $z_{13}$, $z_{23}$ all lie in this strip.
However, by the holomorphicity of the function $r^\tau$,
(\ref{3.37}) is then necessarily valid for any $\omega$, $z$ for which
$r^\tau$ is defined by eqs. (\ref{3.35}), (\ref{3.36}). 

Of course, 
it is possible to calculate  $(\pi_{v_1} \otimes \pi_{v_2})( r^{k,-}(\omega))$ 
as well on an appropriate domain of $v_1$, $v_2$. This is left as an exercise.

\subsection{Recovering Felder's $r$-matrices}

In this subsection $\G$ is a complex simple Lie algebra,
and we start by fixing a Cartan subalgebra and a corresponding set 
$\Phi^+$ of positive roots. 
We also choose root vectors 
$E_\alpha$ ($\alpha\in \Phi$) and dual bases of the Cartan 
subalgebra, $H_i$ and $H^j$, normalized so that
\begin{equation}
B(H_i, H^j)= \delta_i^j,
\qquad
B(E_\alpha, E_{-\alpha})= 1.
\label{3.55}\end{equation}
If $\alpha_i\in \Phi^+$ are the simple roots,
then there is a unique element, $J$, of the Cartan subalgebra 
for which
\begin{equation}
\alpha_i(J)=1
\qquad
\forall i=1,\ldots, \mathrm{rank}(\G).  
\label{3.56}\end{equation}
Let $N$ be the largest eigenvalue of $(\ad J)$ plus $1$,
i.e., the Coxeter number of $\G$.
We wish to show that the application of our preceding 
construction to the automorphism
\begin{equation}
\mu:= \exp(\frac{2\pi i}{N} \ad J)
\label{3.57}\end{equation}
provides an $r$-matrix that is equivalent to Felder's 
solution of the CDYBE \cite{Feld}.
The fixed point set $\G_0$ of this $\mu$ is the chosen Cartan subalgebra of $\G$,
and Felder's $r$-matrix is in fact equivalent to 
\begin{equation}
S^\tau(\omega, z):= \frac{1}{2\pi i} 
\frac{\theta_1'(z\vert \tau)}{\theta_1(z\vert \tau)} H_i \otimes H^i+
\sum_{\alpha\in \Phi} \chi(\alpha(\omega), z\vert \tau)
E_\alpha \otimes E_{-\alpha}.
\label{3.58}\end{equation}
To be precise,  
Felder's original $r$-matrix, ${\cal F}^\tau$,  is given by 
${\cal F}^\tau(\omega,z):= 2\pi i S^\tau(2\pi i \omega, z)$,
which is a substitution that leaves the CDYBE invariant.
Referring to the corresponding terms in (\ref{3.58}), 
below we also write 
$S^\tau:= S^\tau_{\mathrm{Cartan}}+ S^\tau_{\mathrm{root}}$.

It is well known that $\mu$ (\ref{3.57}) acts as a Coxeter element 
on a Cartan subalgebra  which is `in opposition' to the
Cartan subalgebra $\G_0$ and that $\A(\G,\mu)$ 
with its natural gradation is isomorphic to the untwisted
affine Lie algebra of $\G$ equipped with its principal gradation \cite{Kac}. 
In \cite{EV} the homogeneous realization of the untwisted affine Lie
algebra was used to recover Felder's $r$-matrix with 
the aid of evaluation homomorphisms.  
The principal realization provided by $\A(\G,\mu)$ 
must of course give an equivalent result.
It is enlightening to see how this works, 
and it also provides a useful check on our foregoing calculations.

By using the above notations, now we can spell out 
$r^\tau$ from (\ref{3.35}), (\ref{3.36}) explicitly as 
$r^\tau= r^\tau_{\mathrm{Cartan}} + r^\tau_{\mathrm{root}}$ with
\begin{equation}
r^\tau_{\mathrm{Cartan}}(\omega,z)=
B(H_i, \chi_0( \ad \omega, z\vert \tau) H_j) H^i \otimes H^j
= \chi_0(0, z\vert \tau) H_i \otimes H^i.
\label{3.59}\end{equation}
The second equality holds because 
$\chi_0( \ad \omega, z\vert \tau) H_j = \chi_0(0, z\vert \tau) H_j$,
which in turn follows from $(\ad \omega) H_j=0$.
It is easy to compute that
\begin{equation}
\chi_0(0, z\vert \tau)= \lim_{w\rightarrow 0} \chi_0(w,z\vert \tau)=
\frac{1}{2\pi i} 
\frac{\theta_1'(z\vert \tau)}{\theta_1(z\vert \tau)}.
\label{3.60}\end{equation}
Thus the Cartan parts of $S^\tau$ and $r^\tau$ are equal, and
are $\omega$-independent.

As for the root part, by using that
$(\ad \omega) E_\alpha= \alpha(\omega) E_\alpha$,
the definitions give
\begin{eqnarray}
&& r^\tau_{\mathrm{root}}(\omega,z)=
\sum_{\alpha \in \Phi^+} e^{\frac{2\pi i \alpha(J) z}{N}} 
\chi( \alpha(\omega) + 2\pi i \frac{\alpha(J)}{N} \tau , z\vert \tau)
E_\alpha \otimes E_{-\alpha}\qquad\qquad\qquad\qquad
\nonumber\\
&&\phantom{r^\tau_{\mathrm{root}}(\omega,z)} +
\sum_{\alpha \in \Phi^+} e^{\frac{2\pi i (N-\alpha(J)) z}{N}}
\chi( -\alpha(\omega) + 2\pi i \frac{N-\alpha(J)}{N} \tau , z\vert \tau)
E_{-\alpha} \otimes E_{\alpha}.
\label{3.61}\end{eqnarray}
Then we use the identity
\begin{equation}
\chi(w+2\pi i \tau,z\vert \tau)= e^{-2\pi i z} \chi(w,z\vert\tau),
\label{3.62}\end{equation}
which permits us to rewrite $r^\tau_{\mathrm{root}}$ as
\begin{equation}
 r^\tau_{\mathrm{root}}(\omega,z)=
\sum_{\alpha \in \Phi} e^{\frac{2\pi i \alpha(J) z}{N}} 
\chi( \alpha(\omega) + 2\pi i \frac{\alpha(J)}{N} \tau , z\vert \tau)
E_\alpha \otimes E_{-\alpha}.
\label{3.63}\end{equation}

By comparing the above expressions of $r^\tau$ and $S^\tau$, we conclude that
\begin{equation}
r^\tau(\omega,z)= 
\left(e^{ \frac{2\pi i}{N} z_1 \ad J}\otimes  
e^{ \frac{2\pi i}{N} z_2 \ad J}\right)
S^\tau(\omega + 2\pi i \frac{\tau}{N} J, z\vert \tau)
\quad
\hbox{with}\quad z=z_1-z_2.
\label{3.64}\end{equation}
If the dynamical variable $\omega$ belongs to a Cartan subalgebra,
$\G_0$,
then the constant shifts of $\omega$ and the similarity transformations 
by $e^{z_1 \ad H} \otimes e^{z_2 \ad H}$ for any $H\in \G_0$,
$z_1-z_2=z$ map the solutions of the CDYBE to other solutions.
In fact, these transformations are special cases of the gauge
transformations considered in section 4.2 of \cite{EV}.

In summary, we have shown that 
the solution of the CDYBE provided by proposition 2 
in the principal case of $\mu$ in (\ref{3.57}) is gauge equivalent
to Felder's dynamical $r$-matrix in the sense of (\ref{3.64}).

Recently generalizations of Felder's $r$-matrices have been 
found \cite{ES3} for which the dynamical variables 
belong to a subalgebra of a Cartan of a simple Lie algebra $\G$.
The subalgebra in question is the fixed point set of an 
outer automorphism of $\G$ of finite order, and the $r$-matrices given by 
proposition 4.2 in \cite{ES3} contain the same elliptic 
functions that appear in (\ref{3.35}).
These $r$-matrices are very likely to be     
gauge equivalent to those special cases of the $r$-matrices 
constructed in subsection 3.3 for which $\G$ is 
simple and $\G_0$ is a contained in a Cartan subalgebra.  
The precise relationship will be described elsewhere.

\section{Conclusion}
\setcounter{equation}{0}

The purpose of this paper has been to further develop
the construction of dynamical $r$-matrices 
building mainly on the seminal paper \cite{EV} and 
our recent work \cite{FGP1}.
Here our first main result is theorem 1,
whereby a dynamical $r$-matrix  is associated with any 
graded self-dual Lie algebra subject to the rather mild 
conditions in (\ref{2.2})--(\ref{2.4}) and the strong
spectral condition described in (\ref{2.6}).
Our second main result is the application of   
this construction to the general class of affine Lie algebras
$\A(\G,\mu)$ corresponding to the automorphisms
of the finite-dimensional self-dual Lie algebras 
that preserve the scalar product and are of finite order.
The resulting dynamical $r$-matrices are 
generalizations of the basic trigonometric dynamical $r$-matrices of \cite{EV},
which are recovered if $\mu$ is a Coxeter automorphism of a simple Lie algebra.
Motived by the derivation of Felder's elliptic dynamical $r$-matrices \cite{Feld} 
found in \cite{EV},
we have also determined the
spectral-parameter-dependent $\G\otimes \G$-valued dynamical $r$-matrices
that correspond to the $\A(\G,\mu)\otimes\A(\G,\mu)$-valued 
$r$-matrices directly obtained from theorem 1.
The result is given explicitly by proposition 2 and 
proposition 3 is subsection 3.3.
 
It is worth noting that the conditions of theorem 1 are 
satisfied also if $\A$ is an arbitrary Kac-Moody Lie algebra 
associated with a symmetrizable generalized Cartan matrix,
equipped  with the principal gradation \cite{Kac}. In this case 
one recovers the $r$-matrices given by equation (3.4) in \cite{EV}.
It would be interesting to find applications of theorem 1
outside the aforementioned classes of Lie algebras.
As candidates, we plan to examine the two-dimensional (toroidal) 
analogues of the affine Lie algebras \cite{EF}.

Another interesting problem is to find applications 
of the generalizations of Felder's $r$-matrices provided by
proposition 2 in integrable systems.  
In this respect, it appears promising to seek for   
generalized Calogero-Moser type systems, 
since certain spin Calogero-Moser systems 
are known to be closely related to Felder's $r$-matrices \cite{ABB,Xu}.
A different possibility is to uncover these $r$-matrices in the framework of
generalized WZNW models, such as those introduced recently in \cite{Klim}.
This approach would require generalizing  
the results in \cite{BFP} about the exchange $r$-matrices of the usual
WZNW model.
We wish to  pursue this line of research in the future.

\bigskip\bigskip\medskip
\noindent{\bf Acknowledgements.}
We are indebted to J. Balog and A. G\'abor for
reading the manuscript and for discussions.
We are also grateful to P. Etingof for drawing our
attention to reference \cite{ES3}.
This investigation was supported in part by the Hungarian 
Scientific Research Fund (OTKA) under T034170, 
T029802, T030099  and M028418. 
\medskip

\renewcommand{\theequation}{\arabic{section}.\arabic{equation}}
\renewcommand{\thesection}{\Alph{section}}
\setcounter{section}{0} 

\section{The maximal open domain $\check \K\subset \A(\G,\mu)_0$}
\setcounter{equation}{0}
\renewcommand{\theequation}{A.\arabic{equation}}

In this appendix we show that if $\A=\A(\G,\mu)$, then 
the {\em maximal, nonempty, open} domain 
on which the $r$-matrix of theorem 1 can be defined is given by 
$\check \K$ in (\ref{3.11}), where $\omega\in \B_k$ is
subject to the conditions in (\ref{3.12}) and (\ref{3.13}).

In general, the elements of the domain $\check \K \subset \A_0$ 
must satisfy the spectral conditions (\ref{2.6}).
If $\A=\A(\G,\mu)$ and $\kappa \in \K$ is parametrized as 
in (\ref{3.10}), then these conditions are explicitly given by (\ref{3.12}) and
(\ref{3.13}), where $\lambda_a$ is an arbitrary eigenvalue of $\ad \omega\vert \G_a$.
Since $\lambda_0=0$ is always one of the eigenvalues, 
the second condition in (\ref{3.13}) implies that 
$k\neq 2\pi i \frac{n}{m}$ for any $n\in \bZ$, $m\in \bZ^*$. 
As $\check \K$ must be an {\em open} subset of $\K$, 
it follows that $k\in (\bC \setminus i \bR)$ 
for any admissible $\kappa= \omega + kd + l\hat c$.
Note that $\check \K \neq \emptyset$, since e.g.  
the elements of the form 
$\kappa= k d + l\hat c$ in (\ref{3.11}) satisfy the 
conditions (\ref{3.12}), (\ref{3.13}).
Hence we only have to show that (\ref{3.11}) subject to these conditions 
is an {\em open} subset of $\K$.

If $\lambda_a$ is an arbitrary eigenvalue of $\ad \omega$ on $\G_a$ and
$k\in (\bC\setminus i\bR)$, then let us consider the real 
line in $\bC$ defined by 
\begin{equation}
L_{\lambda_a, k}(t)= \lambda_a + kt,
\qquad
\forall t\in \bR.
\label{A.1}\end{equation}
This line intersects the imaginary axis for $t=t_{\lambda_a,k}$ 
at the point $P_{\lambda_a,k}=L_{\lambda_a,k}(t_{\lambda_a,k})$,
\begin{equation}
t_{\lambda_a,k} = -\frac{\Re(\lambda_a)}{\Re(k)},
\qquad
P_{\lambda_a,k} = \lambda_a-k\frac{\Re(\lambda_a)}{\Re(k)}.
\label{A.2}\end{equation}
Now the condition in (\ref{3.12}) can be reformulated as follows: 
\begin{equation}
 P_{\lambda_a,k}\notin 2\pi i \bZ \qquad\hbox{or}\qquad 
 t_{\lambda_a,k} \notin (a + N \bZ),
\qquad\forall a\in \E_\mu \setminus \{0\}.
\label{A.3}\end{equation}
This can be further  reformulated as the requirement 
\begin{equation}
\left\vert e^{P_{\lambda_a,k}} -1 \right\vert^2  +
\left\vert e^{\frac{2\pi i}{N}( t_{\lambda_a,k}-a)} -1 \right\vert^2 \neq 0.
\label{A.4}\end{equation}
It is also useful to rephrase 
the second condition in (\ref{3.13}) as 
\begin{equation}
 P_{\lambda_0,k}\notin 2\pi i \bZ \qquad\hbox{or}\qquad 
 t_{\lambda_0,k} \notin N \bZ^*.
\label{A.5}\end{equation}
Let $\T: \bC \rightarrow \bC$ be an arbitrary continuous function,
which is zero precisely on $N\bZ^*$. 
(For example, we may use $\T(z)=z^{-1}\sin(N^{-1} \pi z)$.)
Then (\ref{A.5}) is equivalent to
 \begin{equation}
\left\vert e^{P_{\lambda_0,k}} -1 \right\vert^2 +
\left\vert \T(t_{\lambda_0,k}) \right\vert^2\neq 0.
\label{A.6}\end{equation}
Since the left hand sides of 
(\ref{A.4}) and (\ref{A.6})  are given by continuous functions of $k$ and
the $\lambda_a$, it follows that these inequalities are stable 
with respect to small variations of $k$ and the $\lambda_a$.
The same is true for the first condition $\lambda_0\notin 2\pi i \bZ^*$
in (\ref{3.13}).
The statement that $\check \K\subset \K$ and $\B_k\subset \G_0$
subject to (\ref{3.11}), (\ref{3.12}), (\ref{3.13}) are {\em open}
subsets follows from this observation by taking into account 
that the position of the eigenvalues of $\ad \omega$ varies 
continuously with $\omega\in \G_0$. 
This means that by choosing $\omega$ near enough to say $\omega^*$,
any eigenvalue of $\ad \omega$ can be taken to be arbitrarily close to
some eigenvalue of $\ad \omega^*$.

\section{A remark on some finite-dimensional $r$-matrices}
\setcounter{equation}{0}
\renewcommand{\theequation}{B.\arabic{equation}}

We here describe some finite-dimensional dynamical
$r$-matrices, which were first considered in the appendix of \cite{ES2},
and point out a relationship between these and the
infinite-dimensional $r$-matrices described in subsection 3.2.

Let $\mu$ be an automorphism of a self-dual Lie algebra
of the same type as in section 3 and recall the
decomposition in (\ref{3.1}), (\ref{3.2}).
For any $a\in \E_\mu$ and integer $q$ specified below,
introduce the meromorphic function $f_{a,q}$ by
\begin{equation}
f_{0,q}(w):= \frac{1}{2}\coth\frac{w}{2}- \frac{1}{w},
\quad
f_{a,q}(w):= \frac{1}{2}\coth\frac{1}{2}(w+\frac{2\pi i}{N} qa)
\quad\quad\hbox{if}\quad 
a\neq 0. 
\label{B.1}\end{equation}
In order to guarantee that these functions are holomorphic
in a neighbourhood of $w=0$, we require the integer
 $q$ to satisfy the conditions
\begin{equation}
1\leq q\leq (N-1), \qquad
q a\notin N \bZ^* 
\quad 
\forall a\in \E_\mu\setminus \{0\}.
\label{B.2}\end{equation} 
Then there exists a nonempty open domain $\check \G_0 \subset \G_0$,
containing the origin,
on which the map $\rho_q: \check \G_0 \rightarrow \mathrm{End}(\G)$
can be defined by
\begin{equation}
\rho_q(\omega) \xi := f_{a,q}(\ad \omega) \xi
\qquad \forall \xi \in \G_a,
\quad
\omega\in \check \G_0.
\label{B.3}\end{equation}

It can be shown that 
$\rho_q$ satisfies the CDYBE (\ref{1.3}),
where $\A$ is replaced by $\G$ 
and $\K$ is taken to be $\G_0$.
If $\mu=\mathrm{id}$, then $\rho_q$ becomes the well known 
canonical (or Alekseev-Meinrenken) dynamical $r$-matrix \cite{EV,AM,BFP}.
In the case $q=1$, which always satisfies (\ref{B.2}), 
$\rho_q$ has been introduced in \cite{ES2}, where it was proved 
that it solves the CDYBE.
The proof given in \cite{ES2} is very elegant
and is very indirect.
A direct proof in the case $\mu = \mathrm{id}$ 
is written down in \cite{PF1}. For general $\mu$ and $q$,
a proof of the CDYBE for $\rho_q$ can be extracted from 
the following observation.
If we let 
$k:=\frac{2\pi i}{N}q$,
then we have 
\begin{equation}
\rho_q(\omega) \eta = R_k(\omega)\eta
\quad\hbox{and}\quad 
(\rho_q(\omega) \xi)^{n_a}= R_k(\omega) \xi^{n_a} 
\label{B.4}\end{equation}
for any $\eta\in \G_0$ and $\xi\in \G_a$, 
$a\neq 0$, $n_a\in (a+ N \bZ)$, 
where $R_k$ refers to the formula (\ref{3.16}).
It should be stressed that this is a relationship purely
at the level of formulas, since in the definition of
the infinite-dimensional $r$-matrices in section 3  
the imaginary values of $k$ were 
excluded for domain reasons.
Nevertheless, it follows from this coincidence 
of formulas that essentially the same algebraic 
computation that proves the CDYBE (\ref{3.17})
can be repeated to verify the CDYBE for $\rho_q$.
We have also verified the CDYBE for $\rho_q$ by
a direct calculation that proceeds analogously to the
proof of our theorem 1.

In certain cases
$\rho_q$ is equivalent to an $r$-matrix of the form in (\ref{1.2})
by a shift of the dynamical variable.
Namely, this happens if the automorphism $\mu$ can be written as 
\begin{equation}
\mu=\exp( \frac{2\pi i}{N} \ad M),
\qquad
M\in \G,
\label{B.5}\end{equation}
where $\ad M$ is diagonalizable and 
the fixed point set $\G_0$ of $\mu$ satisfies 
\begin{equation}  
\G_0= \mathrm{Ker}( \ad M).
\label{B.6}\end{equation}
In particular, by (\ref{B.5}), $\mu$ is an inner automorphism of $\G$.
If these assumptions hold, then we can define a new $r$-matrix 
$\tilde \rho_q$ by
\begin{equation} 
\tilde \rho_q(\omega):= \rho_q(\omega -\frac{2\pi i}{N}q M),
\label{B.7}\end{equation} 
and this $r$-matrix can be identified with $R$ in (\ref{1.2}) 
by taking $\A:=\G$ and $\K:= \G_0$.
The $\G_0$-equivariance property of the dynamical $r$-matrices 
is respected by the shift of the variable in (\ref{B.7})  
on account of (\ref{B.6}).


\newpage
\begin{thebibliography}{99}

\bibitem{EV} P. Etingof and A. Varchenko, 
Commun. Math. Phys. 192 (1998) 77, q-alg/9703040.

\bibitem{Feld}
G. Felder, pp.~1247-1255 in: Proc. Int. Congr. Math. Z\"urich,   
1994, hep-th/9407154;\\
G. Felder and C. Wieczerkowski, 
Commun. Math. Phys. 176 (1996) 133, hep-th/9411004.

\bibitem{ABB}
J. Avan, O. Babelon and E. Billey, 
Commun. Math. Phys.  178 (1996) 281, hep-th/9505091.

\bibitem{Xu}
L.C. Li and P. Xu, 
Integrable spin Calogero-Moser systems,
math.QA/0105162.

\bibitem{ES1}
P. Etingof and O. Schiffmann,
Lectures on the dynamical Yang-Baxter equations,
math.QA/9908064.

\bibitem{Dub}
L. Feh\'er,
Dynamical $r$-matrices and the chiral WZNW phase space,
to appear in the proceedings of the
`Group23' International Colloquium, math-ph/0104027.

\bibitem{Figu}
J.M. Figueroa-O'Farrill and S. Stanciu,
J. Math. Phys. 37 (1996) 4121, hep-th/9506152.

\bibitem{FGP1}
L. Feh\'er, A. G\'abor and B.G. Pusztai,
J. Phys. A  34 (2001) 7235, math-ph/0105047. 

\bibitem{AM} 
A. Alekseev and E. Meinrenken,
Invent. Math. 139 (2000) 135, math.DG/9903052.

\bibitem{BFP}
J. Balog, L. Feh\'er and L. Palla, 
Phys. Lett. B 463 (1999) 83, hep-th/9907050;\\
J. Balog, L. Feh\'er and L. Palla, 
Nucl. Phys. B 568 (2000) 503, hep-th/9910046.

\bibitem{ES2}
P. Etingof and O. Schiffmann, 
Math. Res. Lett.  8 (2001) 157, math.QA/0005282.

\bibitem{DS}
N. Dunford and J.T. Schwartz, Linear operators, Part I: General theory,  
Interscience Publ. Inc., New York-London, 1958.

\bibitem{SW}
D.H. Sattinger and O.L. Weaver, 
Lie groups and algebras with applications to physics, geometry, 
and mechanics, Springer-Verlag, New York-Berlin, 1986.

\bibitem{WW}
E.T. Whittaker and G.N. Wattson,
A course of modern analysis, fourth edition,
Cambridge University Press, Cambridge, 1927. 

\bibitem{Kac}
V.G. Kac, Infinite dimensional Lie algebras, third edition, 
Cambridge University Press, Cambridge, 1990.

\bibitem{ES3}
P. Etingof and O. Schiffmann, 
Math. Res. Lett. 6 (1999) 593, math.QA/9908115.

\bibitem{EF}
P. Etingof and I. Frenkel, Commun. Math. Phys. 165 
(1994) 429, hep-th/9303047.

\bibitem{Klim}
C. Klimcik, Quasitriangular WZW model, hep-th/0103118.

\bibitem{PF1}
B.G. Pusztai and L. Feh\'er,
J. Phys. A  34 (2001) 10949, math.QA/0109082.

\end{thebibliography}
\end{document}